\documentclass[12pt,a4paper]{article}

\usepackage{amssymb}

\def\cM{{\cal M}}
\def\Pic{{\rm Pic}}

\def\oM{{\overline{\cal M}}}

\def\qed{{\hfill $\diamondsuit$}}

\def\Aut{{\rm Aut}}
\def\Hir{{\rm Hir}}
\def\KP{{\rm KP}}
\def\oKP{{\overline{\rm KP}}}
\def\LKP{{\rm LKP}}
\def\Z{{\mathbb Z}}
\def\C{{\mathbb C}}
\def\d{{\partial}}
\def\bF{{\textit{\textbf{F}}\,}}
\def\bL{{\textit{\textbf{L}}\,}}
\def\bl{{\textit{\textbf{l}}\,}}
\def\sbl{{\textit{\textbf{\scriptsize l}}\,}}

\usepackage{theorem}

\newtheorem{theorem}{Theorem}
\newtheorem{proposition}{Proposition}[section]
\newtheorem{corollary}[proposition]{Corollary}
\newtheorem{lemma}[proposition]{Lemma}
\newtheorem{conjecture}[proposition]{Conjecture}
{\theorembodyfont{\rmfamily}
\newtheorem{definition}[proposition]{Definition}
\newtheorem{example}[proposition]{Example}
\newtheorem{remark}[proposition]{Remark}

}

\include{epsf}

\title{Changes of variables in ELSV-type formulas}

\author{Sergey Shadrin\thanks{
Department of Mathematics,
University of Zurich,
Winterthurerstrasse 190,
CH-8057 Zurich, Switzerland
and
Department of Mathematics,
Institute of System Research,
Nakhimovsky prospekt 36-1,
Moscow 117218, Russia.
E-mails:
sergey.shadrin@math.uzh.ch or
shadrin@mccme.ru\,.
Partly supported by the grants RFBR-05-01-01012a, NSh-4719.2006.1,
NWO-RFBR-047.011.2004.026 (RFBR-05-02-89000-NWO-a), 
by the G{\"o}ran Gustafsson foundation, and by
Pierre Deligne's fund based on his 2004 Balzan prize in mathematics.},
Dimitri Zvonkine\thanks{
Institut math{\'e}matique de Jussieu,
Universit{\'e} Paris~VI, 175, rue du Chevaleret,
75013 Paris, France. E-mail: zvonkine@math.jussieu.fr.
Partly supported by the ANR project ``Geometry and
Integrability in Mathematical Physics'' ANR-05-BLAN-0029-01.}}

\date{\today}

\begin{document}

\maketitle

\section{Introduction}
\label{Sec:Intro}

In~\cite{GoJaVa} I.~P.~Goulden, D.~M.~Jackson, and R.~Vakil 
formulated a conjecture relating certain Hurwitz
numbers (enumerating ramified coverings of the sphere) to the
intersection theory on a conjectural Picard variety $\Pic_{g,n}$.
This variety, of complex dimension $4g-3+n$, is supposedly endowed with
a natural morphism to the moduli space of stable curves
$\oM_{g,n}$. The fiber over a point $x \in \cM_{g,n}$ lying
in the open part of the moduli space is equal to the jacobian
of the corresponding smooth curve $C_x$. The variety $\Pic_{g,n}$
is also supposed to carry a universal curve ${\cal C}_{g,n}$ with $n$
disjoint sections $s_1, \dots, s_n$. Denote by ${\cal L}_i$
the pull-back under $s_i$ of the cotangent
line bundle to the fiber of ${\cal C}_{g,n}$. Then we obtain 
$n$ tautological 2-cohomology classes $\psi_i = c_1 ({\cal L}_i)$
on $\Pic_{g,n}$.

We are going to use Goulden, Jackson, and Vakil's 
formula to study the intersection
numbers of the classes $\psi_i$ on $\Pic_{g,n}$
(if it is ever to be constructed). 
In particular, we prove a Witten-Kontsevich-type 
theorem relating the intersection theory
and integrable hierarchies. These equations,
together with the string and dilaton equations, allow us
to compute all the intersection numbers under consideration.

Independently of the conjecture
of~\cite{GoJaVa}, our results can be interpreted as meaningful
statements about Hurwitz numbers.

Our methods are close
to those of M.~Kazarian and S.~Lando
in~\cite{KazLan} and make use of Hurwitz numbers.
We also extend the results of~\cite{KazLan} to include the Hodge
integrals over the moduli spaces, involving one $\lambda$-class.

\subsection{The conjecture}

Fix $n$ positive integers $b_1, \dots, b_n$. Let
$d = \sum b_i$ be their sum. 

\begin{definition} \label{Def:HurwitzPic}
The number of degree~$d$ ramified coverings of the sphere by a genus~$g$
surface possessing a unique preimage of~$0$, $n$ numbered preimages of $\infty$
with multiplicities $b_1, \dots, b_n$, and $2g-1+n$ fixed simple
branch points is called a {\em Hurwitz number} and denoted by
$h_{g;b_1, \dots, b_n}$.
\end{definition}

\begin{conjecture} \label{Conj:GJV}
{\rm (I.~P.~Goulden, D.~M.~Jackson, R.~Vakil, \cite{GoJaVa})}
There exists a compactification of the Picard variety over
the moduli space $\cM_{g,n}$ by a smooth $(4g-3+n)$-dimensional orbifold 
$\Pic_{g,n}$, natural cohomology classes 
$\Lambda_2, \dots, \Lambda_{2g}$ on $\Pic_{g,n}$
of {\rm (}complex{\rm) } degrees $2, \dots, 2g$,
 and an extension of the tautological classes 
$\psi_1, \dots, \psi_n$ such that
$$
h_{g; b_1, \dots, b_n} = (2g-1+n)! \, d \;
\int_{\Pic_{g,n}}
\frac{1 - \Lambda_2 + \cdots \pm \Lambda_{2g}}
{(1-b_1 \psi_1) \cdots (1-b_n \psi_n)} .
$$
\end{conjecture}

Assuming that the conjecture is true we can define
\begin{equation} \label{Eq:bracketPic}
\left< \tau_{d_1} \cdots \tau_{d_n} \right> = 
\int_{\Pic_{g,n}} \psi_1^{d_1} \cdots \psi_n^{d_n}.
\end{equation}
By convention, this bracket vanishes unless $\sum d_i = 4g-3+n$.
We also introduce the following generating series for the intersection
numbers of the $\psi$-classes on $\Pic_{g,n}$:
\begin{equation} \label{Eq:FPic}
F(t_0, t_1, \dots) = \sum_n \frac{1}{n!}
\sum_{d_1, \dots, d_n}
\left< \tau_{d_1} \cdots \tau_{d_n} \right>
t_{d_1} \cdots t_{d_n}.
\end{equation}
We denote by 
\begin{equation} \label{Eq:U}
U = \frac{\d^2 F}{\d^2 t_0} 
\end{equation}
its second partial derivative.

While Conjecture~\ref{Conj:GJV} remains open, the situation should be
seen in the following way. The Hurwitz numbers turn out
to have the unexpected property of being polynomial in
variables~$b_i$ (first conjectured in~\cite{GouJac2} and 
proved in~\cite{GoJaVa}). The coefficients of these polynomials
are denoted by 
$\left< \tau_{d_1} \cdots \tau_{d_n} \Lambda_{2k} \right>$.
(We restrict ourselves to the case $k=0$ with $\Lambda_0=1$.)
The conjectured relation of these coefficients with geometry is
a strong motivation to study them. Our goal is to find out
as much as we can about the values of the bracket in the
combinatorial framework, waiting for their geometrical meaning
to be clarified.

This study was, actually, already initiated in~\cite{GoJaVa}.
In particular, the authors proved that the values of the bracket
satisfy the string and dilation equations:
\begin{eqnarray}
\label{Eq:string}
\frac{\d F}{\d t_0} &=& \sum_{d \geq 1} t_d \frac{\d F}{\d t_{d-1}} 
+ \frac{t_0^2}2,\\ 
\label{Eq:dilaton}
\frac{\d F}{\d t_1} &=& 
\frac12 \sum_{d \geq 0} (d+1) t_d \frac{\d F}{\d t_d}
- \frac12 F .
\end{eqnarray}

By abuse of language we will usually speak of the coefficients
of~$F$ as intersection numbers, implicitly assuming the
conjecture to be true.

\subsection{Results}

We will soon see that $F$ is related to the following
generating function for the Hurwitz
numbers:
\begin{equation} \label{Eq:H}
H(\beta, p_1, p_2, \dots) = \sum_{g,n} \frac{1}{n!}
\frac{\beta^{2g-1+n}}{(2g-1+n)!} 
\sum_{b_1, \dots, b_n}
\frac{h_{g;b_1, \dots, b_n}}{d} p_{b_1} \cdots p_{b_n}.
\end{equation}
Here, as before, $d= \sum b_i$ is the degree of the coverings
and $2g-1+n$ is the number of simple branch points.

Denote by $L_p$ the differential operator 
$$
L_p = \sum b p_b \frac{\d}{\d p_b}.
$$
Its action on $H$ consists in multiplying each term by its total
degree~$d$.

\begin{theorem} \label{Thm:HirotaH}
The series $L_p^2 H$ is a  $\tau$-function
of the Ka\-dom\-tsev--Petviashvili {\rm(}or KP\/{\rm)} hierarchy in
variables~$p_i$; i.e., it satisfies the full set of bilinear Hirota
equations. In addition, $L_p^2 H$ satisfies
the linearized KP equations.
\end{theorem}

The proof of this theorem follows in an almost standard way
from the general theory of integrable systems. We will
discuss it in Section~\ref{Sec:Hirota}.

\begin{theorem} \label{Thm:HirotaU}
The series $U$ is a $\tau$-function
for the KP hierarchy in variables $T_i = t_{i-1}/(i-1)!$\/; i.e., 
it satisfies the full set of bilinear Hirota equations
in these variables. In addition, it satisfies
the linearized KP equations in the same variables.
\end{theorem}

Theorem~\ref{Thm:HirotaU} follows from Theorem~\ref{Thm:HirotaH}, 
but far from trivially, in spite of  their apparent similarity.

\begin{example} 
The string and the dilaton equations allow one
to compute all the values of the bracket in $g=0,1,2$ knowing
only the following values, which can be obtained using
Theorem~\ref{Thm:HirotaU}:
$$
\begin{array}{ll}
g=0: \quad & \displaystyle \left< \tau_0^3 \right> = 1. \\
\\
g=1: & \displaystyle \left< \tau_2 \right> = \frac1{24}. \\
\\
g=2: & \displaystyle \left< \tau_6 \right> = \frac1{1920}, \;\;
       \left< \tau_2 \tau_5 \right> = \frac{19}{5760}, \;\;
       \left< \tau_3 \tau_4 \right> = \frac{11}{1920}, \;\;
       \left< \tau_2^2 \tau_4 \right> = \frac{37}{1440},\\
\\
&      \displaystyle
       \left< \tau_2 \tau_3^2 \right> = \frac5{144}, \;\;
       \left< \tau_2^3 \tau_3 \right> = \frac5{24}, \;\;
       \left< \tau_2^5 \right> = \frac{25}{16}.
\end{array}
$$
\end{example}

\subsection{Acknowledgments}

We are grateful to M.~Kazarian and S.~Lando for stimulating
discussions. We also thank the Stockholm University, where
the major part of this work done, for its hospitality.

\section{Intersection numbers and Hurwitz numbers}
\label{Sec:IntHur}

Here we establish a link between the generating series
$H$ (for Hurwitz numbers) and $F$ (for intersection numbers
of the $\psi$-classes on $\Pic_{g,n}$).

Introduce the following linear triangular change of variables:
\begin{equation} \label{Eq:Chvar}
p_b = \sum_{d=b-1}^{\infty}  \beta^{-(d+1)/2} \,
\frac{(-1)^{d-b+1}}{(d-b+1)! (b-1)!} \; t_d.
\end{equation}
Thus
$$
\begin{array}{lcccccrcl}
p_1 &=& \beta^{-1/2} t_0 &- & \beta^{-1} t_1 & + & \frac12 \beta^{-3/2} t_2
& - & \cdots \; ,\\
p_2 &=&              &    & \beta^{-1} t_1 & - & \beta^{-3/2} t_2 
& - & \cdots \; ,\\
p_3 &=&                     & &        & & \frac12 \beta^{-3/2} t_2 
& - & \cdots \; .\\
\end{array}
$$

Let us separate the generating series $H$ into~2 parts.
The {\em unstable} part, corresponding to the cases $g=0$,
$n=1,2$, is given by
$$
H_{\rm unst}(\beta, p_1, p_2 \dots) =
\sum_{b=1}^{\infty} \frac{p_b}{b^2}
+ \frac{\beta}2 \sum_{b_1,b_2=1}^{\infty} \frac{p_{b_1} p_{b_2}}{b_1+b_2} .
$$
The {\em stable} part is given by $H_{\rm st} = H - H_{\rm unst}$.

The change of variables was designed to make the following proposition
work.

\begin{proposition} \label{Prop:Chvar1}
The change of variables~{\rm (\ref{Eq:Chvar})} transforms the
series $H_{\rm st}$ into a series of the form 
$\sqrt{\beta} \, F + O(\beta)$.
\end{proposition}

\paragraph{Proof.} First let $\beta=1$. It is readily seen that, 
for any $d \geq 0$,
\begin{equation} \label{Eq:binomial}
\sum_{b=1}^{d+1} \frac{(-1)^{d-b+1}}{(d-b+1)! (b-1)!} 
\cdot \frac1{1 - b  \psi}
= \psi^d + O(\psi^{d+1})
\end{equation}
as a power series in~$\psi$. Using Conjecture~\ref{Conj:GJV},
it follows that for any $d_1, \dots, d_n$ we have
$$
\sum_{
\begin{array}{c} 
 \scriptstyle b_1, \dots, b_n \\ 
 \scriptstyle 1 \leq b_i \leq d_i+1
\end{array}
}
\!\!\!\!\!\!\!
\frac{(-1)^{d-b+1}}{(d-b+1)! (b-1)!} \; 
\frac{h_{g;b_1, \dots, b_n}}{(2g-1+n)! \, d}
\qquad \qquad \qquad \qquad
$$
$$
\qquad \qquad \qquad \qquad
= \int\limits_{\Pic_{g,n}} \!\!\!
(1 - \Lambda_2 + \cdots \pm \Lambda_{2g}) \, 
\prod_{i=1}^n (\psi_i^{d_i} + O(\psi_i^{d_i+1})).
$$
Now {\em assume} that $\sum d_i = \dim(\Pic_{g,n}) = 4g-3+n$.
Then each factor in the right-hand side contributes to the
integral only through its lowest order term. Therefore the
right-hand side is equal to
$$
\int\limits_{\Pic_{g,n}} \psi_1^{d_1} \cdots \psi_n^{d_n},
$$
which is, up to a combinatorial factor, precisely a coefficient 
of~$F$. The purpose of introducing the parameter $\beta$ in the
change of variables~(\ref{Eq:Chvar}) 
is precisely to isolate such terms from
the others. Indeed, we claim that the power of $\beta$ in a term
obtained by the change of variables equals 
$$
\frac{\dim(\Pic_{g,n}) - \sum d_i + 1}2.
$$
To check this, recall that the power of $\beta$ in a term
of $H$ equals $2g-1+n$ by definition of $H$. After subtracting
$(d+1)/2$ for each variable $t_d$ we obtain 
$$
2g-1+n/2 - \sum d_i/2 = \frac{4g-2+n - \sum d_i}2 = 
\frac{\dim(\Pic_{g,n}) - \sum d_i + 1}2
$$
as claimed. Thus applying the change of variables to~$H$ we
obtain a series with only positive (half-integer) powers 
of~$\beta$, and the lowest order terms in $\beta$ form the
series $\sqrt{\beta} F$.
\qed

\bigskip

The transformation of the partial derivatives corresponding
to~(\ref{Eq:Chvar}) is obtained by computing the inverse matrix.
It is given by
\begin{equation} \label{Eq:invmatrix}
\frac{\d}{\d p_b} = 
\sum_{d=0}^{b-1} \beta^{(d+1)/2} \frac{(b-1)!}{(b-d-1)!} \; \frac{\d}{\d t_d}.
\end{equation}
Thus
$$
\begin{array}{rcrcrcrc}
\displaystyle \frac{\d}{\d p_1} 
&=& 
\displaystyle \beta^{1/2} \frac{\d}{\d t_0},\\
\\
\displaystyle \frac{\d}{\d p_2} 
&=& 
\displaystyle \beta^{1/2} \frac{\d}{\d t_0}
&+&
\displaystyle \beta \, \frac{\d}{\d t_1},\\
\\
\displaystyle \frac{\d}{\d p_3} 
&=& 
\displaystyle \beta^{1/2} \frac{\d}{\d t_0}
&+&
\displaystyle 2 \beta \, \frac{\d}{\d t_1}
&+&
\displaystyle 2 \beta^{3/2} \frac{\d}{\d t_2}.\\
\end{array}
$$

\begin{proposition} \label{Prop:Chvar2}
The change of variables~{\rm (\ref{Eq:Chvar})} induces the following
transformations:
$$
\begin{array}{rcl}
L_p & \longrightarrow & \displaystyle
L_t = \sum\limits_{d \geq 0} (d+1) t_d \frac{\d}{\d t_d} 
+ \frac1{\sqrt{\beta}} \sum\limits_{d \geq 1} t_d \frac{\d}{\d t_{d-1}} \;,\\
\\
L_p^2 H_{\rm unst} & \longrightarrow & \displaystyle
\frac1{\sqrt{\beta}} \, t_0 (t_1+1) + t_0^2.
\end{array}
$$
\end{proposition}

Both claims of the proposition are obtained by simple computations. \qed

\medskip

The concinnity of this result is striking. Indeed, both transforms 
could have contained arbitrarily large negative powers of
$\beta$, but they happen to cancel out in both cases.
Further, $L_p^2 H_{\rm unst}$ is an infinite series, but
after the change of variables it has become
a polynomial with only three terms. 
Most important of all, the coefficients $L_{-1}$ and
$L_0$ of $\beta^{-1/2}$ and $\beta^0$, respectively, 
in the operator $L_t$ are precisely the string and 
dilaton operators from Equations~(\ref{Eq:string}) 
and~(\ref{Eq:dilaton}). This leads to the following corollaries.

\begin{corollary} \label{Cor:1}
We have
$$
L_t F = F + 2 \frac{\d F}{\d t_1} + 
\frac1{\sqrt{\beta}} \left(\frac{\d F}{\d t_0} -\frac{t_0^2}2 \right).
$$
$$
L_t \frac{\d F}{\d t_0} = 2 \frac{\d^2 F}{\d t_0 \d t_1} + 
\frac1{\sqrt{\beta}} \left(\frac{\d^2 F}{\d t_0^2} -t_0 \right).
$$
\end{corollary}

\begin{corollary} \label{Cor:HtoU}
The change of variables~{\rm (\ref{Eq:Chvar})} induces the
following transformations of generating series:
$$
\begin{array}{rcl}
L_p^2 H_{\rm st} & \longrightarrow & \frac1{\sqrt{\beta}} \, 
\Bigl[U - t_0(t_1+1) \Bigr] + O_\beta(1) \;,\\
\\
L_p^2 H & \longrightarrow & \frac1{\sqrt{\beta}} U + O_\beta(1) \;.\\
\end{array}
$$
Here $O_\beta(1)$ is a series containing only nonnegative powers 
of~$\beta$.
\end{corollary}

\paragraph{Proof.} The first result follows from Corollary~\ref{Cor:1},
while the second one is obtained after a (yet another!) cancellation
of the term $t_0(t_1+1)$ with the contribution of 
$L_p^2 H_{\rm unst}$. \qed

\section{Hirota equations and KP hierarchy}
\label{Sec:Hirota}
\label{Ssec:HirWedge}

In this section we recall some necessary facts about the
Hirota and the KP hierarchies and use them to
prove Theorem~\ref{Thm:HirotaH}.

We start with a brief introduction to the Hirota and the
KP hierarchy. More details can be found, for instance,
in~\cite{KacRai}.

The {\em semi-infinite wedge space} $W$ is the vector space
of formal (possibly infinite) linear combinations of
infinite wedge products of the form
$$
z^{k_1} \wedge z^{k_2} \wedge \dots,
$$
with $k_i \in \Z$, $k_i = i$ starting from some~$i$. Consider a sequence
$\varphi_1, \varphi_2, \dots$ of Laurent series 
$\varphi_i \in \C[[z^{-1},z]$ such that 
$\varphi_i = z^i + \mbox{(lower order terms)}$ starting
from some~$i$. Then $\varphi_1 \wedge \varphi_2 \wedge \dots$
is an element of~$W$. The elements
that can be represented in that way are called {\em decomposable}.
One way to check whether an element is
decomposable, is to verify if it satisfies the Pl{\"u}cker
equations. 

Now we will assign an element of $W$ to any power series in 
variables $p_1, p_2, \dots$. The series will turn out to
be a solution of the Hirota hierarchy if and only if the corresponding
element of~$W$ is decomposable.

To a Young diagram $\mu$ with $d$ squares we assign the
{\em Schur polynomial} $s_\mu$ in variables $p_1, p_2, \dots$
defined by
$$
s_\mu = \frac1{d!} \sum_{\sigma \in S_d} \chi_\mu (\sigma) p_\sigma. 
$$
Here $S_d$ is the symmetric group, $\sigma$ is a permutation,
$\chi_\mu (\sigma)$ is the character of $\sigma$ in the
irreducible representation assigned to $\mu$, and
$p_\sigma = p_{l_1} \cdots p_{l_k}$, where $l_1, \dots, l_k$
are the lengths of cycles of $\sigma$. 

The Schur polynomials $s_\mu$ with $\mbox{area}(\mu)=d$
form a basis of the space of quasihomogeneous polynomials
of weight~$d$ (the weight of $p_i$ being equal to~$i$).

Consider a power series $\tau$ in variables $p_i$. Decomposing it
in the basis of Schur polynomials we can uniquely assign to it
a (possibly infinite) linear combination of Young diagrams $\mu$
(of all areas). Now, in this linear combination 
we replace each Young diagram $\mu = (\mu_1, \mu_2, \dots)$ 
by the following wedge product:
$$
z^{1-\mu_1} \wedge z^{2-\mu_2} \wedge \dots,
$$
where $(\mu_1, \mu_2, \dots)$ are the lengths of the columns
of $\mu$ in decreasing order with an infinite number of zeroes
added in the end. We have obtained an element $w_\tau \in W$.

The bilinear Pl{\"u}cker equations on the coordinates of $w_\tau$
happen to combine into bilinear differential equations on~$\tau$,
called the Hirota equations. Thus, as we said, $w_\tau$ is
decomposable if and only if $\tau$ is a solution of the Hirota
equations. Let us define these equations precisely.

Consider two partitions $\lambda$ and $\mu$ of an integer~$d$.
Denote by $\chi_\mu (\lambda)$ the character of any
permutation with cycle type $\lambda$ in the irreducible
representation assigned to~$\mu$.
Denote by $|\Aut(\lambda)|$ the number of permutations of
the elements of $\lambda$ that preserve their values. For instance
$|\Aut(7,6,6,4,1,1,1,1,1)| = 2! \cdot 5!$.

Let $d_i = \d/\d p_i$.
Let $D_\mu$ be the differential operator
$$
D_\mu = \sum_{\lambda, \; |\lambda|=d} \chi_\mu (\lambda)
\frac{d_{\lambda_1} \cdots d_{\lambda_k}}{|\Aut (\lambda)|},
$$
where $k$ is the number of elements 
of~$\lambda$. For instance
$$
D_{()} = 1, \quad D_{(1)} = d_1, \quad D_{(2)} = \frac12 d_1^2 + d_2,
\quad D_{(1,1)} = \frac12 d_1^2 - d_2,
$$
$$
D_{(3)} = \frac16 d_1^3 + d_1 d_2 + d_3, \quad 
D_{(2,1)} = \frac13 d_1^3 - d_3 , \quad 
D_{(1,1,1)} = \frac16 d_1^3 - d_1 d_2 + d_3.
$$

Let $\tau$ be a formal power series in $p_1, p_2, \dots$. If
the constant term of $\tau$ does not vanish, we can 
also consider its logarithm $F = \ln \tau$.

\begin{definition} (see~\cite{DubNat}, Proposition~1).
\label{Def:hierarchies}
The {\em Hirota hierarchy} is the following family of 
bilinear differential equations:
\begin{equation}
\Hir_{i,j}(\tau) = D_{()}\tau \cdot D_{(j,i)}\tau 
- D_{(i-1)}\tau \cdot D_{(j,1)}\tau
+ D_{(j)}\tau \cdot D_{(i-1,1)}\tau,
\end{equation}
for $2 \leq i \leq j$.

Substituting $\tau = e^F$ and dividing by $\tau^2$
we obtain a family of equations
on~$F$. It is called the {\em Kadomtsev--Petviashvili} or
{\em KP} hierarchy:
\begin{equation}
\KP_{i,j}(F) = \frac{\Hir_{i,j}(e^F)}{e^{2F}}.
\end{equation}

Finally, leaving only the linear terms in the KP hierarchy
we obtain the {\em linearized KP equations}\footnote{In the
published version of this paper, the linearized KP equations are
erronously called the {\em dispersionless limit} of KP.}:
\begin{equation}
\LKP_{i,j}(F) = \mbox{linear part of } \KP_{i,j}(F).
\end{equation}
\end{definition}

\begin{example}
Denoting the derivative with respect to $p_i$ by the
index~$i$, we have
\begin{eqnarray*}
\Hir_{2,2} &=& 
\tau \tau_{2,2} -\tau_2^2- \tau \tau_{1,3} + \tau_1 \tau_3 +
\frac14 \tau_{1,1}^2 - \frac13 \tau_1 \tau_{1,1,1} +
\frac1{12} \tau \tau_{1,1,1,1},\\
\KP_{2,2}&=& F_{2,2} - F_{1,3} + \frac12 F_{1,1}^2 + 
\frac1{12} F_{1,1,1,1},\\
\LKP_{2,2}&=& F_{2,2} - F_{1,3} + \frac1{12} F_{1,1,1,1},\\
\Hir_{2,3} &=& \tau \tau_{2,3} - \tau_2 \tau_3 - \tau \tau_{1,4}
+ \tau_1 \tau_4 +\frac12 \tau_{1,1} \tau_{1,2} 
- \frac12 \tau_1 \tau_{1,1,2} - \frac16 \tau_{1,1,1} \tau_2 \\
&&+ \frac16 \tau \tau_{1,1,1,2} - \frac12 \tau_1 \tau_{2,2}
+ \frac12 \tau \tau_{1,2,2} - \tau_{1,2} \tau_2
- \frac12 \tau \tau_{1,1,3} + \frac12 \tau_{1,1} \tau_3\\
&& + \frac1{24} \tau \tau_{1,1,1,1,1} - \frac18 \tau_1 \tau_{1,1,1,1}
+ \frac1{12} \tau_{1,1} \tau_{1,1,1},\\
\KP_{2,3} &=& F_{2,3} - F_{1,4} + F_{1,1} F_{1,2} + \frac16
F_{1,1,1,2} + \frac12 F_1 F_{2,2} - \frac12 F_1 F_{1,3}
+ \frac14 F_1 F_{1,1}^2 \\
&& + \frac1{24} F_1 F_{1,1,1,1}
+ \frac12 F_{1,2,2} - \frac12 F_{1,1,3} + \frac12 F_{1,1} F_{1,1,1} + 
\frac1{24} F_{1,1,1,1,1},\\
\LKP_{2,3} &=& F_{2,3} - F_{1,4}  + \frac16
F_{1,1,1,2} + \frac12 F_{1,2,2} - \frac12 F_{1,1,3}  
+ \frac1{24} F_{1,1,1,1,1}.\\
\end{eqnarray*}
\end{example}

\begin{remark}
Every Hirota equation can be simplified by adding to it
some partial derivatives of lower equations. This, in turn, leads to
simplified KP and LKP equations. For instance, we have
$$
\begin{array}{lcl}
\displaystyle \Hir_{2,3} -\frac12 \frac{\d \Hir_{2,2}}{\d p_1} &=& 
\displaystyle \tau \tau_{2,3} - \tau_2 \tau_3 - \tau \tau_{1,4}
+ \tau_1 \tau_4 + \frac12 \tau_{1,1} \tau_{1,2} \\
&&\displaystyle
- \frac12 \tau_1 \tau_{1,1,2} - \frac16 \tau_{1,1,1} \tau_2
+ \frac16 \tau \tau_{1,1,1,2},\\
\\
\displaystyle 
\KP_{2,3} - \frac12 F_1 \; \KP_{2,2} - \frac12 \frac{\d \KP_{2,2}}{\d p_1}
&=& \displaystyle 
F_{2,3} - F_{1,4} + F_{1,1} F_{1,2} + \frac16 F_{1,1,1,2},\\
\\
\displaystyle 
\LKP_{2,3} - \frac12 \frac{\d \LKP_{2,2}}{\d p_1} &=&
\displaystyle 
F_{2,3} - F_{1,4}  + \frac16 F_{1,1,1,2}.
\end{array}
$$
Thus we obtain a simplified hierarchy which is, of course,
equivalent to the initial one.
Sometimes it is the equations of this simplified 
hierarchy that are called Hirota equations.
However, for our purposes it is easier to use 
the equations as we defined them.
\end{remark}

\paragraph{Proof of Theorem~\ref{Thm:HirotaH}.} We will
actually prove that {\em for any function $c= c(\beta)$,
the series $c + L_p^2 H$ satisfies the Hirota hierarchy}.
Let us show that the element of $W$ assigned to
$c+L_p^2 H$ is decomposable.
Consider the following Laurent series in~$z$:
$$
\varphi_1 = cz+\sum_{n \geq 0} \beta^{n(n+1)/2} z^{-n},
\qquad 
\varphi_i = z^i - e^{(i-1) \beta} z^{i-1}
\quad \mbox{for } i \geq 2. 
$$
The coefficients of these series are shown in the matrix below.
$$
\begin{array}{cccccccccccc}
\dots & -4 & -3 & -2 & -1 & 0 & 1 & 2 & 3 & 4 & 5 &\dots\\
\hline
\\
\dots & e^{10 \beta} & e^{6 \beta} & e^{3 \beta} & e^{\beta} &
1 & c & 0 & 0 & 0 & 0 & \dots \\
\dots & 0 & 0 & 0 & 0 &
0 & -e^{- \beta} & 1 & 0 & 0 & 0 & \dots \\
\dots & 0 & 0 & 0 & 0 &
0 & 0 & -e^{-2\beta} & 1 & 0 & 0 & \dots \\
\dots & 0 & 0 & 0 & 0 &
0 & 0 & 0 & -e^{-3 \beta} & 1 & 0 & \dots \\
\dots & 0 & 0 & 0 & 0 &
0 & 0 & 0 & 0 & - e^{-4 \beta} & 1 & \dots
\end{array}
$$
We claim that expanding the wedge product 
$\varphi_1 \wedge \varphi_2 \wedge \dots$
and replacing every Young diagram by the corresponding
Schur polynomial we obtain the series $c+ L_p^2 H$. The
proof goes as in~\cite{KazLan}.

We introduce the so-called {\em cut-and-join operator}
$$
A = \frac12 \sum_{i,j=1}^\infty
\left[
(i+j) p_i p_j \frac{\d}{\d p_{i+j}}
+ ij p_{i+j} \frac{\d^2}{\d p_i \d p_j}
\right] \; .
$$
Then $L_pH$ satisfies the equation $\d (L_p H)/ \d \beta = A (L_pH)$
(see~\cite{GouJac}).
Since the operator $L_p$ commutes both with $A$ and 
with $\d/ \d \beta$, the series $L_p^2 H$ satisfies
the same equation. 

The Schur polynomials $s_\lambda$ are eigenvectors of~$A$. The eigenvalue
corresponding to a Young diagram $\lambda$ equals
$f_\lambda = \frac12 \sum \lambda_i (\lambda_i - 2i +1)$, where $\lambda_i$
are the column lengths.

This allows one to reconstitute the
whole series $L_p^2 H$ starting with its $\beta$-free terms
$L_p^2 H|_{\beta=0}$: if
$$
L_p^2 H|_{\beta=0} = \sum c_\lambda s_\lambda
$$
then 
$$
L_p^2 H = \sum c_\lambda s_\lambda e^{f_\lambda \beta}.
$$

It is apparent from the form of the above matrix
that the coefficients of $s_\lambda$ in the expansion
are nonzero only in two cases: (i)~for the empty diagram,
where the coefficient equals $c$; (ii)~for the {\em hook Young diagrams} 
$\lambda = \mbox{hook}(a,b)$ with column lengths 
$$
a+1, \underbrace{1, 1, \dots, 1}_b.
$$
For a Young diagram like that, the coefficient of
$s_{\mbox{\scriptsize hook}(a,b)}$ equals
$$
(-1)^b \, e^{[a(a+1)/2 - b(b+1)/2]\beta}.
$$
For the $\beta$-free terms we have
$$
L_p^2 H |_{\beta=0} = \sum_{i \geq 1} p_i
= \sum_{a,b \geq 0} (-1)^b s_{\mbox{\scriptsize hook}(a,b)},
$$
the second equality being an exercise in the representation theory.
To this we add the remark that $a(a+1)/2 - b(b+1)/2$ is
precisely the eigenvalue $f_\lambda$ for 
$\lambda = \mbox{hook}(a,b)$. It follows 
that the series corresponding to 
$\varphi_1 \wedge \varphi_2 \wedge \dots$
equals $L_p^2 H$ as claimed.

Thus we have proved that the series $c + L_p^2 H$ 
satisfies the Hirota hierarchy.

The claim about the $\LKP$ equations is a simple corollary of that.
Indeed, it follows from Definition~\ref{Def:hierarchies}
that for any series $G$ we have
$$
\Hir_{i,j}(1+G) - \Hir_{i,j}(G) = \LKP_{i,j}(G).
$$
Thus, from the fact that both $L_p^2H$ and $1+L_p^2H$ satisfy the
Hirota equations it follows immediately that $L_p^2H$
satisfies the linearized KP equations. 
\qed

\section{Hirota equations and the change of variables}
\label{Sec:HirChvar}

Now we are going to study the effect of the change of
variables~(\ref{Eq:Chvar}) on the Hirota hierarchy and 
prove Theorem~\ref{Thm:HirotaU}.

In Section~\ref{Ssec:HirWedge} we assigned to each Young
diagram $\mu$ an operator $D_\mu$ and used these operators
as building blocks to define the Hirota equations. It turns
out that the change of variables~(\ref{Eq:Chvar}) acts
on $D_\mu$ by ``biting off'' the corners of $\mu$. From this
we will deduce that each Hirota equation becomes, after the
change of variables, a linear combination of lower Hirota
equations.

\bigskip

As in Theorem~\ref{Thm:HirotaU}, we rescale the variable
$t_i$ by setting $t_i = i! T_{i+1}$. Then we have
\begin{equation} \label{Eq:T}
\frac{\d}{\d p_i} = {i-1 \choose 0} \beta^{1/2} \frac{\d}{\d T_1}
+ {i-1 \choose 1} \beta \frac{\d}{\d T_2} + \cdots
+{i-1 \choose i-1} \beta^{i/2} \frac{\d}{\d T_i}
\end{equation}
instead of~Eq.~(\ref{Eq:invmatrix}).

Consider a linear differential operator $D$ with constant coefficients
in variables $p_i$. Denote $\d/\d p_i$ by $d_i$ and consider $d_i$
as a new set of variables. Introduce the differential
operator
$$
S = \sum_{i \geq 1} i d_i \frac{\d}{\d d_{i+1}}
$$
in these variables. Then we have the following lemma.

\begin{lemma} \label{Lem:S}
Assume $D$ is a quasi-homogeneous polynomial in variables $d_i$
with total weight~$n$. Then applying the change of 
variables~{\rm (\ref{Eq:T})} to~$D$ viewed as a differential operator
is equivalent
to applying the differential operator $\beta^{n/2}e^{S/\sqrt{\beta}}$
to $D$ viewed as a polynomial in variables~$d_i$.
\end{lemma}

The proof is a simple check. \qed

\medskip

What we actually want is to apply the change of variables~(\ref{Eq:T})
to the operators $D_\mu$ defined in Section~\ref{Ssec:HirWedge}.
Indeed, these operators are the building blocks of the Hirota
equations (Definition~\ref{Def:hierarchies}). The answer is
given below in Proposition~\ref{Prop:SDmu}.

\begin{definition} \label{Def:corner}
A square of a Young diagram $\mu$ is called a {\em corner} if, when
we erase it, we obtain another Young diagram. In other words,
a corner is a square with coordinates
$(i,\mu_i)$ such that either $\mu_{i+1} < \mu_i$ or $\mu_i$ is the
last column of~$\mu$.
If $(i, \mu_i)$ is a corner of $\mu$ we will denote by
$\mu - \square_i$ the diagram obtained by erasing this corner.
\end{definition}

\begin{proposition} \label{Prop:SDmu}
We have 
$$
S D_\mu = \!\!\!\!\!\!\sum_{(i,\mu_i) = \mbox{\small \rm corner of } \mu}
\!\!\!\!\!\!\!
(\mu_i - i) \cdot  D_{\mu - \square_i} \; .
$$
\end{proposition}

\paragraph{Proof.}
Let $\mu$ be a Young diagram with $d$ squares
and $\lambda$ a partition of~$d-1$. Assume that $\lambda$
has $k$ parts. Then, for $1 \leq i \leq k$, denote by
$\lambda+1_i$ the partition of $d$ obtained from $\lambda$
by replacing $\lambda_i$ by $\lambda_i+1$.

Writing down explicitly the action of
$S$ on $D_\mu$ one finds that the assertion of the proposition
is equivalent to the following identity:
\begin{equation} \label{Eq:identity}
\sum_{(i,\mu_i) = \mbox{\small \rm corner of } \mu}
\!\!\!\!\!\!\!
(\mu_i -i) \cdot \chi_{\mu-\square_i}(\lambda)
=
\sum_{i=1}^k \lambda_i \cdot \chi_\mu (\lambda+1_i) \; .
\end{equation}

To prove this identity we need a short digression into the
representation theory of the symmetric group as presented
in~\cite{OkoVer}.

Consider the subgroup $S_{d-1} \subset S_d$ consisting of the
permutations that fix the last element~$d$.
The irreducible representation of $S_d$ assigned 
to~$\mu$ is then also a representation of $S_{d-1}$, although
not necessarily irreducible. It turns out that this representation
is isomorphic to 
$$
\bigoplus_{\mbox{\small corners of } \mu} \mu-\square_i,
$$
where, by abuse of notation, $\mu-\square_i$ stands for the
irreducible representation of $S_{d-1}$ assigned to this Young
diagram.

Further, consider the following element of the group
algebra $\C S_d$:
$$
X = (1,d) + (2,d) + \cdots + (d-1,d).
$$
This element (the sum of all transpositions involving $d$)
is called the first Jucys--Murphey-Young element. It obviously
commutes with the subgroup $S_{d-1}$. Therefore its eigenspaces
in the representation $\mu$
coincide with the irreducible subrepresentations of $S_{d-1}$,
{\em i.e.}, they are also in one-to-one correspondence with
the corners of~$\mu$. The eigenvalue corresponding to the 
corner $(i,\mu_i)$ equals $\mu_i-i$ (see~\cite{OkoVer}).

Using this information, let us choose a permutation $\sigma \in S_{d-1}$
with cycle type $\lambda$ and compute in two different ways
the character
$$
\chi_\mu (\sigma \cdot X),
$$
where $\sigma \cdot X \in \C S_d$.

First way. Both $\sigma$ and $X$ leave invariant the irreducible
subrepresentations of $S_{d-1}$. For $X$ such a
subrepresentation is an eigenspace with eigenvalue
$\mu_i-i$. The character of $\sigma$ in the same subrepresentation
equals $\chi_{\mu - \square_i} (\lambda)$. We obtain the
left-hand side of Identity~(\ref{Eq:identity}).

Second way. Let us see what happens when we multiply $\sigma$
by~$X$. Each transposition in $X$ increases the length of
precisely one cycle of $\sigma$ by~$1$. This is equivalent to
increasing one of the $\lambda_i$'s by~1. Moreover if the $i$th
cycle of $\sigma$ has length $\lambda_i$, it will be touched by
a transposition from~$X$ exactly $\lambda_i$ times. Thus we obtain
the right-hand side of Identity~(\ref{Eq:identity}).

This completes the proof. \qed

\begin{proposition} \label{Prop:SHir}
The change of variables~{\rm (\ref{Eq:T})} transforms the Hirota
equation $\Hir_{i,j}$ into an equation of the form
$$
\sum_{
\begin{array}{c}
\scriptstyle 2 \leq i' \leq i, \;\; 2 \leq j' \leq j\\
\scriptstyle i' \leq j'
\end{array}
}
c_{i',j'} \; \beta^{(i'+j')/2} \; \Hir_{i', j'} \; 
$$
for some rational constants $c_{i', j'}$. The constant $c_{i,j}$
of the leading term equals~$1$.
\end{proposition}

\paragraph{Proof.} By Definition~\ref{Def:hierarchies} the equation
$\Hir_{i,j}$ has the form
$$
\Hir_{i,j}(\tau) = D_{()}\tau \cdot D_{(j,i)}\tau 
- D_{(i-1)}\tau \cdot D_{(j,1)}\tau
+ D_{(j)}\tau \cdot D_{(i-1,1)}\tau.
$$
According to Lemma~\ref{Lem:S}, applying the change 
of variables to the equation is the same as applying to 
each $D_\mu$ in this expression the operator
$$
\beta^{(i+j)/2}e^{S/\sqrt{\beta}}.
$$

To simplify the computations, consider the flow $e^{tS}$ applied
to $\Hir_{i,j}$. We will compute the derivative of this flow
with respect to~$t$. If $\cal D$ is the vector space of all
polynomials in variables $d_i$, then $\Hir_{i,j}$ lies in 
$\cal D \otimes \cal D$. The flow $e^{tS}$ acts as
$e^{tS} \otimes e^{tS}$, while its derivative with respect to~$t$
is $1 \otimes S + S \otimes 1$.

We will prove that $1 \otimes S + S \otimes 1$ applied to
$\Hir_{i,j}$ is a linear combination
of lower Hirota equations ($i'< i, j' < j$). 
Since this is true for all $i,j$,
when we integrate the flow we see that $\Hir_{i,j}$ will have 
changed by a linear combination of lower Hirota equations.

It remains to apply $1 \otimes S + S \otimes 1$ to $\Hir_{i,j}$.
To do that we use Proposition~\ref{Prop:SDmu}.

If $i<j$ we obtain
$$
\begin{array}{rcl}
D_{()} \cdot 
\biggl[(i-2) D_{(j,i-1)} + (j-1) D_{(j-1,i)}\biggr] &
+ & 0  \cdot  D_{(j,i)}\\ 
- \; D_{(i-1)}  \cdot (j-1)D_{(j-1,1)} &
- & (i-2) D_{(i-2)}  \cdot  D_{(j,1)}\\
+ \; D_{(j)} \cdot (i-2)D_{(i-2,1)} &
+ & (j-1) D_{(j-1)}  \cdot  D_{(i-1,1)}
\end{array}
$$
$$
= (i-2) \; \Hir_{i-1,j} + (j-1) \;\Hir_{i,j-1}.
$$

If $i=j$ we obtain
$$
\begin{array}{rcl}
D_{()} \cdot (i-2) D_{(i,i-1)}
& + & 0 \cdot D_{(i,i)}\\
- \; D_{(i-1)} \cdot  (i-1) D_{(i-1,1)}
& - & (i-2) D_{(i-2)} \cdot D_{(i,1)}\\
+ \;  D_{(i)}  (i-2) \cdot D_{(i-2,1)}
& + & (i-1) D_{(i-1)} \cdot  D_{(i-1,1)}
\end{array}
$$
$$
= (i-2) \; \Hir_{i-1,i}.
$$

This completes the proof. \qed

\begin{remark} \label{Rem:SHir}
The family of equations given in Proposition~\ref{Prop:SHir}
is equivalent to the Hirota hierarchy. Indeed, the equations
of the Hirota hierarchy can be obtained from these equations
by linear combinations and vice versa.
\end{remark}

\paragraph{Proof of Theorem~\ref{Thm:HirotaU}.}
The series $c+L_p^2 H$ satisfies the Hirota
equations by Theorem~\ref{Thm:HirotaH}. Therefore, by
Proposition~\ref{Prop:SHir} and Remark~\ref{Rem:SHir},
the series obtained from it under the change of 
variables~(\ref{Eq:Chvar}) also satisfies the Hirota hierarchy.
According to Corollary~\ref{Cor:HtoU}, this new series
has the form 
$$
c + \frac{1}{\sqrt{\beta}} U + O_\beta(1).
$$
Taking $c = c'/\sqrt{\beta}$ and considering the lowest
order terms in $\beta$ we obtain that
$c'+U$ satisfies the Hirota hierarchy for any constant~$c'$.
It follows that $U$ satisfies the linearized KP hierarchy. \qed

\section*{Appendix: On Hodge integrals}

\renewcommand{\thesection}{A}
\setcounter{section}{1}

In this section we will use the change of variables suggested
in~\cite{KazLan} to study Hodge integrals over the moduli
spaces of curves. We consider the integrals involving a unique
$\lambda$-class and arbitrary powers of $\psi$-classes.

\subsection{Hurwitz numbers and Hodge integrals}

Here we study intersection theory on moduli spaces
rather than on Picard varieties. We follow the
same path as in Sections~\ref{Sec:Intro} and~\ref{Sec:IntHur},
but with different intersection numbers and Hurwitz numbers.
Our aim is to extend the results of~\cite{KazLan}.

Instead of~(\ref{Eq:bracketPic}) we define the following brackets
\begin{equation} \label{Eq:bracketM}
\left< \tau_{d_1} \cdots \tau_{d_n} \right>^{(k)} = 
\int_{\oM_{g,n}} \psi_1^{d_1} \cdots \psi_n^{d_n} \lambda_k
\end{equation}
for $k+ \sum d_i = 3g-3+n$ (otherwise the bracket vanishes).

Instead of~(\ref{Eq:FPic}) we will use the generating series
\begin{equation} \label{Eq:FM}
F^{(k)}(t_0, t_1, \dots) = \sum_n \frac{1}{n!}
\sum_{d_1, \dots, d_n}
\left< \tau_{d_1} \cdots \tau_{d_n} \right>^{(k)}
t_{d_1} \cdots t_{d_n}.
\end{equation}
We can also regroup these series into a unique series
$$
\bF(z; t_0, t_1, \dots) =  \sum_{k \geq 0} (-1)^k z^k F^{(k)}. 
$$

Instead of the Hurwitz numbers of Definition~\ref{Def:HurwitzPic} 
we now use different Hurwitz numbers.

Fix $n$ positive integers $b_1, \dots, b_n$. Let
$d = \sum b_i$ be their sum. 

\begin{definition} \label{Def:HurwitzM}
The number of degree~$d$ ramified coverings of the sphere by a genus~$g$
surface possessing $n$ numbered preimages of $\infty$
with multiplicities $b_1, \dots, b_n$, and $d+n+2g-2$ fixed simple
branch points is called a {\em Hurwitz number} and denoted by
$h_{g;b_1, \dots, b_n}$.
\end{definition}

We introduce the following generating series for the these numbers:
$$
H(\beta, p_1, p_2, \dots ) = 
\sum_{g,n} \frac{1}{n!} \,
\frac{\beta^{d+n+2g-2}}{(d+n+2g-2)!} 
\sum_{b_1, \dots, b_n}
h_{g;b_1, \dots, b_n} \;\,  p_{b_1} \cdots p_{b_n}.
$$

It is divided in two parts. The unstable part, corresponding to
$g=0$, $n=1,2$ equals
$$
H_{\rm unst} = \sum_{b \geq 1} \beta^{b-1} \frac{b^{b-2}}{b!} p_b
+ \frac12 \sum_{b_1, b_2 \geq 1}
\beta^{b_1+b_2} \frac{b_1^{b_1}b_2^{b_2}}{(b_1+b_2) b_1! b_2 !}
p_{b_1} p_{b_2}.
$$
The stable part equals $H_{\rm st} = H - H_{\rm unst}$.

Finally, instead of Conjecture~\ref{Conj:GJV} we use the so-called
ELSV formula proved in~\cite{ELSV}.

\begin{theorem} {\rm (The ELSV formula~\cite{ELSV})}
We have
$$
h_{g;b_1, \dots, b_n} \; = \; (d+n+2g-2)! \;
\prod_{i=1}^n \frac{b_i^{b_i}}{b_i!} \; 
\int_{\oM_{g,n}} 
\frac{1 - \lambda_1 + \lambda_2 - \cdots \pm \lambda_g}
{(1-b_1 \psi_1) \cdots (1 - b_n \psi_n)} \; .
$$
\end{theorem}

As before, it turns out that the series $\bF$ and $H$ are related
via a change of variables based on Equation~(\ref{Eq:binomial}). 
However the change of variables is different
from~(\ref{Eq:Chvar}) due to (i)~the factors $b_i^{b_i}/b_i!$
in the ELSV formula and (ii)~to a different relation between
the number of simple branch points and the dimension
of the Picard/moduli space. Namely, following~\cite{KazLan}, we let
\begin{equation} \label{Eq:Chvar2}
p_b = \sum_{d=b-1}^\infty
\frac{(-1)^{d-b+1}}{(d-b+1)! \, b^{b-1}} \, \beta^{-b - (2d+1)/3} \; t_d.
\end{equation}
Thus
$$
\begin{array}{lcrcrcrcl}
p_1 &=& \beta^{-4/3} t_0 &- & \beta^{-6/3} t_1 & + & 
\displaystyle \frac12 \beta^{-8/3} \; \, t_2
& - & \cdots \; ,\\
\\
p_2 &=&  & & \displaystyle \frac12 \beta^{-9/3} t_1 & - & 
\displaystyle \frac12 \beta^{-11/3} t_2 
& + & \cdots \; ,\\
\\
p_3 &=& &&& & \displaystyle\frac19 \beta^{-14/3} t_2 
& - & \cdots \; .\\
\end{array}
$$
This change of variables transforms $H$ into a series 
in variables $t_0, t_1, \dots,$ and~$\beta^{2/3}$.

We will also need a more detailed version of Equation~(\ref{Eq:binomial}):
$$
\sum_{b=1}^{d+1} \frac{(-1)^{d-b+1}}{(d-b+1)! (b-1)!} 
\cdot \frac1{1 - b  \psi}
= \psi^d + \sum_{k=1}^\infty a_{d,d+k} \, \psi^{d+k},
$$
where $a_{d,d+k}$ are some rational constants that actually
happen to be integers. For instance,
$$
\begin{array}{lcl}
\displaystyle
\;\;\; \frac1{1-\psi} &=& 1 + \psi + \psi^2 + \cdots,\\
\\
\displaystyle
-\frac1{1-\psi} + \frac1{1-2\psi} 
&=&  \psi + 3\psi^2 + 7 \psi^3 + \cdots,\\
\\
\displaystyle
\;\;\; \frac{1/2}{1-\psi} - \frac1{1-2\psi} + \frac{1/2}{1-3\psi}
&=&  \psi^2 + 6\psi^3 + 25 \psi^4 + \cdots .\\
\end{array}
$$

Using these constants we introduce the
following differential operators:
\begin{eqnarray*}
L_1 &=& \sum_{n=0}^\infty a_{n,n+1} \; t_{n+1} \frac{\d}{\d t_n}\; ,\\
L_2 &=& \sum_{n=0}^\infty a_{n,n+2} \; t_{n+2} \frac{\d}{\d t_n}
+ \frac1{2!}
\sum_{{n_1}, {n_2} = 0}^\infty 
a_{{n_1},{n_1}+1} \, a_{{n_2},{n_2}+1} \; t_{{n_1}+1} t_{{n_2}+1} \; 
\frac{\d^2}{\d t_{n_1} \d t_{n_2}} \; ,\\
L_3 &=& \sum_{n=0}^\infty a_{n,n+3} \; t_{n+3} \frac{\d}{\d t_n}
+ \frac1{2!}
\sum_{{n_1}, {n_2} = 0}^\infty 
a_{{n_1},{n_1}+2} \, a_{{n_2},{n_2}+1} \; t_{{n_1}+2} t_{{n_2}+1} \; 
\frac{\d^2}{\d t_{n_1} \d t_{n_2}} \\
&& + \frac1{2!}
\sum_{{n_1}, {n_2} = 0}^\infty 
a_{{n_1},{n_1}+1} \, a_{{n_2},{n_2}+2} \; t_{{n_1}+1} t_{{n_2}+2} \; 
\frac{\d^2}{\d t_{n_1} \d t_{n_2}} \\
&& + \frac1{3!}
\sum_{{n_1}, {n_2},{n_3} = 0}^\infty \!\!\!\!\!\!
a_{{n_1},{n_1}+1} \, a_{{n_2},{n_2}+1} \, a_{{n_3},{n_3}+1}
\; t_{{n_1}+1} t_{{n_2}+1} t_{{n_3}+1} 
\frac{\d^3}{\d t_{n_1}\d t_{n_2}\d t_{n_3}}\; ,
\end{eqnarray*}
and so on. We can also regroup these operators in a unique operator
$\bL = 1 + zL_1 + z^2 L_2 + \cdots$.

These operators and the change of variables~(\ref{Eq:Chvar2}) were
designed to make the following proposition work.

\begin{proposition} \label{Prop:FandH}
Performing the change of variables~(\ref{Eq:Chvar2}) on the series
$H_{\rm st}$
and replacing $\beta^{2/3}$ by $z$ we obtain the series $\bL \bF$.
\end{proposition}

\paragraph{Proof.} Using the ELSV formula one can check that
the change of variables~(\ref{Eq:Chvar2}) transforms $H$ into
the series
$$
\sum_{n,g} \frac1{n!} \sum_{d_1, \dots d_n}
\!\!\!\! \int\limits_{\ \;\;\;\oM_{g,n}} \!\!\!\!
(1 - \beta^{2/3} \lambda_1 + \beta^{4/3} \lambda_2 - \cdots)
\prod_{i=1}^n
(\psi_1^{d_i} + \beta^{2/3} a_{d_i, d_i+1} \, \psi_1^{d_i+1} + \cdots)
\; .
$$
The proposition follows. \qed

\subsection{Hierarchies and operators}

\begin{proposition} \label{Prop:exp}
We have $\bL = e^\sbl$,
$$
\bl = z l_1 + z^2 l_2 + \cdots
$$
is a first order linear differential operator:
$$
l_k = \alpha_{n,n+k} \; t_{n+k} \frac{\d}{\d t_n}.
$$
\end{proposition}

\paragraph{Proof.}
Consider the operators $l_k$ and $\bl$ as above with indeterminate
coefficients $\alpha_{n,n+k}$.
Consider the expansion of $e^\sbl$ and denote
by $a_{n,n+k}$ the coefficient of $t_{n+k} \, \d/\d t_n$
in this expansion. We have
\begin{eqnarray*}
a_{n,n+1} &=& \alpha_{n,n+1} \; ,\\
a_{n,n+2} &=& \alpha_{n,n+2} + \frac12 \, \alpha_{n,n+1} \,
\alpha_{n+1,n+2} \; ,\\
a_{n,n+3} &=& \alpha_{n,n+3} + \frac12 \, \alpha_{n,n+1} \, 
\alpha_{n+1,n+3}  + \frac12 \, \alpha_{n,n+2} \, \alpha_{n+2,n+3}\\
&& + \frac16 \, \alpha_{n,n+1} \, \alpha_{n+1,n+2} \, \alpha_{n+2,n+3}
\; ,
\end{eqnarray*}
and so on. Note that these equalities allow one to determine the
coefficients $\alpha$ unambiguously knowing the coefficients~$a$.

Now consider the coefficient of a monomial
$$
\prod_{i=1}^p t_{n_i+k_i} \frac{\d}{\d t_{n_i}}
$$
in the same expansion of $e^\sbl$. It is equal to
$$
\frac{|\Aut \{(n_1, k_1), \dots, (n_p, k_p) \}|}{p!}
\prod \alpha_{n_i,n_i+k_i} + \mbox{higher order terms},
$$
and we claim that this sum can be factorized as
$$
\frac{|\Aut \{(n_1, k_1), \dots, (n_p, k_p) \}|}{p!}
\prod a_{n_i,n_i+k_i}.
$$

Indeed, suppose that we have already chosen a power of $\bl$,
say $\bl^q$ and a term in each of the $q$ factors that contribute
to the coefficient of 
$$
\prod_{i=1}^{p-1} t_{n_i+k_i} \frac{\d}{\d t_{n_i}}.
$$
Now we must choose some additional power $\bl^r$ of $\bl$
and a term in each of the $r$ factors that will contribute to
the coefficient of
$$
t_{n_p+k_p} \frac{\d}{\d t_{n_p}} = t_{n+k} \frac{\d}{\d t_n}.
$$
Moreover, we must choose the positions of the $r$ new factors
among the $q$ that are already chosen. This can be done in
$$
{q+r \choose q}
$$
ways. (The operator $t_{n+k} \; \d/\d t_n$
acts by replacing $t_n$ by $t_{n+k}$. Thus the $r$ terms
in question divide the segment $[n,n+k]$ into $r$ parts
and should be ordered in a uniquely determined way.)

In the end we must divide the
coefficient thus obtained by $(q+r)!$ since we are looking at
$e^\sbl$. Thus we obtain a coefficient of 
$$
\frac1{q!} \cdot \frac1{r!}
$$
for any choice of $r$ terms.

Now, if $q=0$ what we have finally obtained is precisely the
expression for $a_{n,n+k}$. For a general $q$ we will therefore
obtain the same expression for $a_{n,n+k}$ divided by $q!$.
Thus we have proved that $a_{n,n+k} = a_{n_p, n_p+k_p}$ 
can be factored out in the coefficient of
$$
\prod_{i=1}^p t_{n_i+k_i} \frac{\d}{\d t_{n_i}}.
$$
The same is true for $a_{n_i, n_i+k_i}$ for all~$i$. Thus
the coefficient is the product of $a_{n_i, n_i+k_i}$ as
claimed.

In other words, we showed that the coefficients of
$\exp (\bl)$ coincide with those of $\bL$. \qed

\begin{conjecture}
The operators $l_k$ have the form
$$
l_k = c_k \sum_{n \geq 0} {n+k+1 \choose k+1} \; t_n 
\frac{\d}{\d t_{n+k}}
$$
for some sequence of rational constants~$c_k$.
\end{conjecture}

The sequence $c_k$ seems quite irregular and starts as follows:
$$
1, \; -\frac12, \; \frac12, \; -\frac23, \; \frac{11}{12}, \; -\frac34, 
\; -\frac{11}6, \; \frac{29}4, \; \frac{493}{12}, \; -\frac{2711}{6},
\; -\frac{12406}{15}, \; \frac{2636317}{60}, \; \dots \, .
$$

\bigskip

Now we will establish a hierarchy of partial differential equations
satisfied by~$\bF$. We use Propositions~\ref{Prop:FandH} 
and~\ref{Prop:exp} together with the following fact.

\begin{theorem} \label{Thm:KP}  {\rm (See~\cite{Okounkov,KazLan})}
The series $H$ satisfies the KP hierarchy in variables
$p_1, p_2, \dots$.
\end{theorem}

This theorem is proved like Theorem~\ref{Thm:HirotaU}
by noticing that $\exp(H)$ satisfies the cut-and-join
equation.

Applying the change of variables~(\ref{Eq:Chvar2}) to
the KP equations we will obtain partial differential
equations satisfied by $\bL \bF$. These equations
can, of course, be considered as equations on $\bF$,
since the coefficients of $\bL$ are known. However,
the equations thus obtained are infinite,
i.e., with an infinite number of terms. Our goal is
to prove that we can combine them in a way that
leads to finite differential equations.

The derivatives $\d /\d p_b$ are expressed via 
$\d / \d t_d$ by computing the inverse of the 
matrix of the change of variables~(\ref{Eq:Chvar2}).
We have
\begin{equation} \label{Eq:deriv}
\frac{\d}{\d p_b} = b^{b-1} \sum_{d=0}^{b-1}
\frac{\beta^{b+(2d+1)/3}}{(b-d-1)!} \, \frac{\d}{\d t_d} \; .
\end{equation}
Thus
$$
\begin{array}{rcrcrcrc}
\displaystyle \frac{\d}{\d p_1} 
&=& 
\displaystyle \beta^{4/3} \frac{\d}{\d t_0},\\
\\
\displaystyle \frac{\d}{\d p_2} 
&=& 
\displaystyle 2 \beta^{7/3} \frac{\d}{\d t_0}
&+&
\displaystyle 2 \beta^{9/3} \frac{\d}{\d t_1},\\
\\
\displaystyle \frac{\d}{\d p_3} 
&=& 
\displaystyle \frac92 \beta^{10/3} \frac{\d}{\d t_0}
&+&
\displaystyle 9 \beta^{12/3} \frac{\d}{\d t_1}
&+&
\displaystyle 9 \beta^{14/3} \frac{\d}{\d t_2}.\\
\end{array}
$$

\bigskip

Now using Theorem~\ref{Thm:KP} and Proposition~\ref{Prop:FandH}
we will transform the KP hierarchy into a system of equations on $\bF$.
We will illustrate the procedure on the example of $\KP_{2,2}$.

We know that $\KP_{i,j}(H) = 0$. For $i=j=2$ this means
$$
\frac{\d^2 H}{\d p_2^2} - \frac{\d^2 H}{\d p_1 \d p_3}
+\frac12 \left( \frac{\d^2 H}{\d p_1^2} \right)^2
+\frac1{12} \frac{\d^4 H}{\d p_1^4} =0.
$$
Using $H = H_{\rm st} + H_{\rm unst}$ and the explicit
expression of $H_{\rm unst}$ we transform $\KP_{i,j}$
into a (finite) equation $\widehat{\rm KP}_{i,j}$ on $H_{\rm st}$.
For instance, for $i=j=2$, we obtain
$$
\frac{\d^2 H_{\rm st}}{\d p_2^2} - \frac{\d^2 H_{\rm st}}{\d p_1 \d p_3}
+\frac12 \left( \frac{\d^2 H_{\rm st}}{\d p_1^2} \right)^2
+\frac1{12} \frac{\d^4 H_{\rm st}}{\d p_1^4} + 
\frac12 \beta^2 \frac{\d^2 H_{\rm st}}{\d p_1^2} = 0  .
$$
Applying the change of variables~(\ref{Eq:deriv}) and replacing
$\beta^{2/3}$ by $z$ we transform this into an equation 
$\oKP_{i,j}$ on $\bL \bF$. For $i=j=2$ we have
$$
- \frac{\d^2 (\bL \bF)}{\d t_0 \d t_1}
+ \frac12 \left( \frac{\d^2 (\bL \bF)}{\d t_0^2} \right)^2
+\frac1{12} \frac{\d^4 (\bL \bF)}{\d t_0^4}
+ z
\left( 
4 \frac{\d^2 (\bL \bF)}{\d t_1^2} - 9 \frac{\d^2 (\bL \bF)}{\d t_0 \d t_2}
\right) = 0 .
$$
In principle, we could have stopped here. However, in this form
the equation is only useful to study the $z$-free part
$F^{(0)}$ of $\bF$, which was done in~\cite{KazLan}. Indeed
the operators $L_i$ for $i \geq 1$ are composed of infinitely
many terms. This means that if we develop the above equation
and take its coefficient of $z^1$, we will obtain an {\em
infinite} equation on $F^{(0)}$ and $F^{(1)}$. Such an equation
is quite useless if we want to compute $F^{(1)}$.
Therefore we continue with the following theorem 
(recall that $\bL = e^\bl$):

\begin{theorem} \label{Thm:finite}
Consider the expression
$$
e^{-\sbl} \;\; \oKP_{i,j}( e^\sbl \bF)
$$
as a series in~$z$. Then its coefficient of
$z^k$ is a finite differential equation on
$F^{(0)}, \dots, F^{(k)}$.
\end{theorem}

\begin{example}
The coefficient of $z^1$ in $e^{-\sbl} \, \oKP_{2,2}( e^\sbl \bF)$
gives the following equation:
$$
- \frac{\d^2 F^{(1)}}{\d t_0 \d t_1}
+ \frac{\d^2 F^{(0)}}{\d t_0^2} \; \frac{\d^2 F^{(1)}}{\d t_0^2}
+ \frac1{12} \frac{\d^4 F^{(1)}}{\d t_0^4} 
\quad \quad \quad \quad \quad \quad
\quad \quad \quad \quad \quad \quad \quad \quad \quad
$$
$$
\quad \quad \quad \quad \quad \quad
\quad \quad \quad 
+12 \frac{\d^2 F^{(0)}}{\d t_0 \d t_2}
-3 \frac{\d^2 F^{(0)}}{\d t_1^2}
-2 \frac{\d^2 F^{(0)}}{\d t_0^2} \; \frac{\d^2 F^{(0)}}{\d t_0 \d t_1}
-\frac13 \frac{\d^4 F^{(0)}}{\d t_0^3 \d t_1} = 0.
$$
Assuming that we know $F^{(0)}$, this equation, together with the
string and the dilaton equations, allow us to compute
all the coefficients of $F^{(1)}$, that is, all Hodge
integrals involving $\lambda_1$.
\end{example}

\paragraph{Proof of Theorem~\ref{Thm:finite}.}
Let $Q$ be a linear differential operator (in variables
$t_d$) whose coefficients are polynomials in~$z$. Then
$$
e^{-\sbl} Q e^\sbl = Q + [Q,\bl] + \frac12[[Q,\bl],\bl] + \cdots
$$
is a series in~$z$ whose coefficients are finite differential
operators. We will denote this series by $\widehat Q$. Now 
suppose we have several linear operators $Q_1, \dots, Q_r$
as above. Since $\bl$ is a first order operator, we obtain
$$
e^{-\sbl} \;\;Q_1(e^\sbl \bF) \cdot \cdots \cdot Q_r(e^\sbl \bF)
= {\widehat Q}_1(\bF) \cdot \cdots \cdot {\widehat Q}_r (\bF).
$$
This is, once again, a series in $z$ whose coefficients are
finite differential equations on the $F^{(k)}$. The
theorem now follows from the fact that every equation
$\oKP_{i,j}$ is a finite linear combination of
expressions of the form
$$
Q_1(e^\sbl \bF) \cdot \cdots \cdot Q_r(e^\sbl \bF).
$$
\qed

Thus every equation $\KP_{i,j}$ and every power
of $z$ gives us a finite differential equation on
the functions $F^{(k)}$. Below we describe some facts 
concerning these equations that we have observed 
but not proved. For any $F^{(k)}$ and any $(d',d'') \not= (0,0)$
by taking linear combinations of the 
equations in question we can obtain an equation
of the form 
$$
\frac{\d F^{(k)}}{\d t_{d'} \d t_{d''}}
= \mbox{terms with more than 2 derivations}.
$$
For homogeneity reasons the sum of indices in these derivatives
will be smaller than $d'+d''$. Therefore we use
similar equations with smaller $d'+d''$ to simplify the
right-hand part by substitutions. After a finite number
of substitutions we will obtain an expression of
$\d F^{(k)}/\d t_{d'} \d t_{d''}$ exclusively via partial
derivatives with respect to $t_0$. Moreover, these expressions
themselves be organized into equations on $\bF$:

\begin{eqnarray*}
\bF_{0,1} &=& 
\left(
\frac12 \bF_{0,0}^2 + \frac1{12}\bF_{0,0,0,0}
\right)
- z
\left(
\frac1{24} \bF_{0,0,0}^2
+ \frac1{720} \bF_{0,0,0,0,0,0}
\right)\\
&& +z^2
\left(
\frac1{720} \bF_{0,0,0} \bF_{0,0,0,0,0}
+\frac1{360} \bF_{0,0,0,0}^2
+\frac1{30240} \bF_{0,0,0,0,0,0,0,0}
\right)\\
&& + \cdots \; , \\
\bF_{0,2} &=& \
\left(
\frac16 \bF_{0,0}^3 
+ \frac1{12} \bF_{0,0} \bF_{0,0,0,0}
+\frac1{24} \bF_{0,0,0}^2 + 
\frac1{240} \bF_{0,0,0,0,0,0}
\right)\\
&&-z
\left(
\frac1{24} \bF_{0,0} \bF_{0,0,0}^2
+\frac1{720} \bF_{0,0} \bF_{0,0,0,0,0,0}
+\frac7{720} \bF_{0,0,0} \bF_{0,0,0,0,0} \right. \\
&& \left. +\frac1{180} \bF_{0,0,0,0}^2
+\frac1{7560} \bF_{0,0,0,0,0,0,0,0}
\right)
+\cdots \; ,\\
\bF_{1,1} &=&  
\left(
\frac13 \bF_{0,0}^3 
+ \frac16 \bF_{0,0} \bF_{0,0,0,0}
+\frac1{24} \bF_{0,0,0}^2 + 
\frac1{144} \bF_{0,0,0,0,0,0}
\right)\\
&& -z
\left(
\frac1{12} \bF_{0,0} \bF_{0,0,0}^2
+\frac1{360} \bF_{0,0} \bF_{0,0,0,0,0,0}
+\frac{13}{720} \bF_{0,0,0} \bF_{0,0,0,0,0} \right. \\
&& \left. +\frac1{120} \bF_{0,0,0,0}^2
+\frac1{4320} \bF_{0,0,0,0,0,0,0,0}
\right)
+\cdots \; , \\
\bF_{0,3} &=& 
\left(
\frac1{24} \bF_{0,0}^4
+\frac1{24} \bF_{0,0}^2 \bF_{0,0,0,0}
+\frac1{24} \bF_{0,0} \bF_{0,0,0}^2
+\frac1{240} \bF_{0,0} \bF_{0,0,0,0,0,0}
\right. \\
&& \left.
+\frac1{120} \bF_{0,0,0} \bF_{0,0,0,0,0} 
+\frac1{160} \bF_{0,0,0,0}^2
+\frac1{6720} \bF_{0,0,0,0,0,0,0,0}
\right)
+ \cdots \; , \\
\bF_{1,2} &=& 
\left(
\frac18 \bF_{0,0}^4
+\frac18 \bF_{0,0}^2 \bF_{0,0,0,0}
+\frac1{12} \bF_{0,0} \bF_{0,0,0}^2
+\frac1{90} \bF_{0,0} \bF_{0,0,0,0,0,0}
\right. \\
&& \left.
+\frac1{60} \bF_{0,0,0} \bF_{0,0,0,0,0} 
+\frac{23}{1440} \bF_{0,0,0,0}^2
+\frac1{2880} \bF_{0,0,0,0,0,0,0,0}
\right)
+ \cdots \; .
\end{eqnarray*}

Not thoroughly unexpectedly, it turns out that the free terms of
these equations form the well-known Korteweg -- de~Vries
(or KdV) hierarchy on $F^{(0)}$.

The first equation above (expressing $\bF_{0,1}$), together
with the string and dilaton equations, is sufficient to 
determine the values of all Hodge integrals involving 
a single $\lambda$-class. This approach seems to be 
simpler than the method of~\cite{Kim} based on the
study of double Hurwitz numbers.

\end{document}